\documentclass[12pt,a4paper]{article}
\usepackage{indentfirst}
\usepackage{makeidx}
\usepackage{amsmath,amssymb,amsfonts,amsthm,graphics,graphicx}
\usepackage{color}
\usepackage{ifpdf}
\usepackage{subfigure}

\newtheorem{thm}{Theorem}

\newtheorem{prob}{Problem}
\newtheorem{lem}{Lemma}
\newtheorem{pro}{Proposition}
\newtheorem{cor}{Corollary}

\theoremstyle{definition}

\def\-{\mbox{--}}

\newtheorem{remark}{Remark}
\newtheorem{obser}{Observation}
\newtheorem{conj}{Conjecture}

\def\pf{\noindent {\it Proof.} }
\begin{document}

\title{\Large\bf Conflict-free vertex-connections of graphs\footnote{Supported by National Science Foundation of China (Nos.
11371205, 11531011, 11601254, 11551001), the Science Found of Qinghai Province
(Nos. 2016-ZJ-948Q, 2014-ZJ-907), and the project on the key lab of
IOT of Qinghai province (No. 2017-Z-Y21).}}
\author{\small Xueliang Li$^{1,2}$, \ \ Yingying Zhang$^{1}$,  \ \  Xiaoyu Zhu$^{1}$,  \ \  Yaping Mao$^{2}$,
\ \ Haixing Zhao$^{2}$ \\[0.2cm]
\small $^{1}$Center for Combinatorics and LPMC\\
\small Nankai University, Tianjin 300071, China\\[0.2cm]
\small $^{2}$School of Mathematics and Statistics\\
\small Qinghai Normal University, Xining, Qinghai 810008, China\\[0.2cm]
\small E-mail: lxl@nankai.edu.cn; zyydlwyx@163.com;\\
\small ~~~~~~~~~zhuxy@mail.nankai.edu.cn; maoyaping@ymail.com; h.x.zhao@163.com}
\date{}
\maketitle
\begin{abstract}
A path in a vertex-colored graph is called \emph{conflict free}
if there is a color used on exactly one of its vertices. A
vertex-colored graph is said to be \emph{conflict-free vertex-connected} if any two vertices of the graph are connected by a conflict-free path. This paper investigates the question: For a connected graph $G$, what is the smallest number of colors needed in a vertex-coloring of $G$ in order to make $G$ conflict-free vertex-connected. As a result, we get that the answer is easy for $2$-connected graphs, and very difficult for connected graphs with more cut-vertices, including trees.


\end{abstract}

\section{Introduction}

In this paper, all graphs considered are simple, finite and undirected. We refer to a book \cite{B} for undefined notation and terminology in graph theory. A path in an edge-colored graph is a {\it rainbow path} if its edges have different colors. An edge-colored graph is {\it rainbow connected} if any two vertices of the graph are connected by a rainbow path of the graph. For a connected graph $G$, the {\it rainbow connection number} of $G$, denoted by $rc(G)$, is defined as the smallest number of colors required to make $G$ rainbow connected. This concept was first introduced by Chartrand et al. in \cite{CJM, CJMZ}. Since then, a lot of
results on the rainbow connection have been obtained; see \cite{LSS, LSu}.

As a natural counterpart of the concept of rainbow connection, the concept of rainbow vertex connection was first introduced by Krivelevich and Yuster in \cite{KY}. A path in a vertex-colored graph is a {\it vertex-rainbow path} if its internal vertices have different colors. A vertex-colored graph is {\it rainbow vertex-connected} if any two vertices of the graph are connected by a vertex-rainbow path of the
graph. For a connected graph $G$, the {\it rainbow vertex-connection number} of $G$, denoted by $rvc(G)$, is defined as the smallest number of colors required to make $G$ rainbow vertex-connected. There are many results on this topic, we refer to \cite{CLS2,LL,LMS1,LS,LMS}.

In \cite{CJV}, Czap et al. introduced the concept of conflict-free
connection. A path in an edge-colored graph is called \emph{conflict free} if there is a color used on exactly one of its edges. An
edge-colored graph is said to be \emph{conflict-free connected} if any two vertices of the graph are connected by a conflict-free path.
The {\it conflict-free connection number} of a connected graph $G$, denoted by $cfc(G)$, is defined as the smallest number of colors required to make $G$ conflict-free connected. Note that for a nontrivial connected graph $G$ with order $n$, we have $$1\leq cfc(G)\leq rc(G)\leq n-1.$$
Moreover, $cfc(G)=1$ if and only if $G$ is a complete graph, and
$cfc(G)=n-1$ if and only if $G$ is a star.

Motivated by the above mentioned concepts, as a natural counterpart
of conflict-free connection number, in this paper we introduce the concept of conflict-free vertex-connection number. A path in a vertex-colored graph is called \emph{conflict free} if there is a color used on exactly one of its vertices. A vertex-colored graph is said to be \emph{conflict-free vertex-connected} if any two vertices of the graph are connected by a conflict-free path. The {\it conflict-free vertex-connection number} of a connected graph $G$, denoted by $vcfc(G)$, is defined as the smallest number of colors required to make $G$
conflict-free vertex-connected. Note that for a nontrivial connected
graph $G$ with order $n$, we can easily observe that
$$2\leq vcfc(G)\leq n.$$
The lower bound is trivial since there is a
path of order at least two between any two vertices in $G$, while
the upper bound is also trivial since one may color all the vertices of
$G$ with distinct colors. The main problem studied in this paper is
the following.

\begin{prob}\label{prob1}
For a given graph $G$, determine its conflict-free vertex-connection number.
\end{prob}

The rest of this paper is organized as follows. In Section $2$, we
prove some preliminary results. In Section $3$, we study the
structure of graphs having conflict-free vertex-connection number
two and three respectively. In Section $4$, we obtain some sharp
bounds of the conflict-free vertex-connection number for trees.

\section{Preliminaries}

The following observation is immediate.
\begin{obser}\label{obs1}
If $G$ is a nontrivial connected graph and $H$ is a connected
spanning subgraph of $G$, then $vcfc(G)\leq vcfc(H)$. In particular,
$vcfc(G)\leq vcfc(T)$ for every spanning tree $T$ of $G$.
\end{obser}

\begin{lem}\label{lem1} Let $G$ be a $2$-connected graph and $w$ be a vertex of $G$. Then for any two vertices $u$ and $v$ in $G$, there is a $u$-$v$ path containing the vertex $w$.
\end{lem}

\pf It is clearly true for the case that $w\in\{u,v\}$ since $G$ is
$2$-connected. Now suppose that $w\in V(G)\backslash\{u,v\}$. Let
$P_1$ and $P_2$ be two internally disjoint paths from $u$ to $w$ in
$G$. If there is a $v$-$w$ path $P$ such that $P$ and $P_1$ are
vertex-disjoint except for the vertex $w$, then the path $uP_1wPv$
is the desired path. Otherwise, let $x$ be the first common vertex
of $P$ and $P_1$ when going along $P$ from $v$. Then the path
$uP_2wP_1xPv$ is the desired path.\qed

For a path, we have the following result.
\begin{thm}\label{thm1}
Let $P_n$ be a path of order $n$. Then $vcfc(P_n)=\lceil\log_2(n+1)\rceil$.
\end{thm}

\pf The proof goes similarly to that of Theorem 3 in \cite{CJV}. Let $P_n=v_1v_2...v_n$. First we show that $vcfc(P_n)\leq\lceil\log_2(n+1)\rceil$. Define a vertex-coloring of $P_n$ by coloring the vertex $v_i$ with color $x+1$, where $i\in[n]$ and $2^x$ is the largest power of $2$ that divides $i$. Clearly, the largest number in such a coloring is $\lceil\log_2(n+1)\rceil$. Moreover, it is easy to check that the maximum color of the vertices on each subpath $Q$ of $P_n$ appears only once on $Q$. Then $P_n$ is conflict-free vertex-connected,  and so $vcfc(P_n)\leq\lceil\log_2(n+1)\rceil$.

Next we just need to prove that $vcfc(P_n)\geq\lceil\log_2(n+1)\rceil$. To show it, it suffices to show that any path with conflict-free vertex-connection number $k$ has at most $2^k-1$ vertices. We apply induction on $k$. The statement is evidently true for $k=2$. Give a path $P_n$ a vertex-coloring with $vcfc(P_n)=k$. Then there is a vertex, say $v_i$, in $P_n$ with a unique color. Delete the vertex $v_i$ from $P_n$. The resulting paths are $P_{i-1}=v_1v_2...v_{i-1}$ and $P_{n-i}=v_{i+1}v_{i+2}...v_n$ with $vcfc(P_{i-1})\leq k-1$ and $vcfc(P_{n-i})\leq k-1$. By the induction hypothesis, $P_{i-1}$ and $P_{n-i}$ have at most $2^{k-1}-1$ vertices, respectively. Thus $P_n$ has at most $2(2^{k-1}-1)+1=2^k-1$ vertices, and so $vcfc(P_n)\geq\lceil\log_2(n+1)\rceil$.

Therefore, $vcfc(P_n)=\lceil\log_2(n+1)\rceil$. \qed

\begin{remark} From Theorem \ref{thm1} and \cite[Theorem 3]{CJV}, we have that $vcfc(P_n)>cfc(P_n)$. However, $vcfc(G)\leq cfc(G)$ if $G$ is a star of order at least $3$. Thus, one of $vcfc(G)$ and $cfc(G)$ cannot be bounded in terms of the other.
\end{remark}

\section{Graphs with conflict-free vertex-connection number two or three}

A {\it block} of a graph $G$ is a maximal connected subgraph of $G$
that has no cut-vertex. Then the block is either a cut-edge, say
{\it trivial block}, or a maximal $2$-connected subgraph, say {\it
nontrivial block}. Let $B_1,B_2,...,B_k$ be the blocks of $G$. The
{\it block graph} of $G$, denoted by $B(G)$, has vertex-set
$\{B_1,B_2,...,B_k\}$ and $B_iB_j$ is an edge if and only if the
blocks $B_i$ and $B_j$ have a cut-vertex in common, where $1\leq
i,j\leq k$.

The following lemma is a preparation of Theorem \ref{thm2}.
\begin{lem}\label{lem2}
Let $G$ be a $2$-connected graph. Then $vcfc(G)=2$.
\end{lem}

\pf Since $vcfc(G)\geq2$, we just need to show that $vcfc(G)\leq2$.
Let $w$ be a vertex of $G$. Define a $2$-coloring $c$ of the
vertices of $G$ by coloring the vertex $w$ with color $2$ and all
the other vertices of $G$ with color $1$. By Lemma \ref{lem1}, for
any two vertices $u$ and $v$ in $G$, there is a $u$-$v$ path
containing the vertex $w$. According to the coloring $c$ of $G$,
this $u$-$v$ path is a conflict-free path. Thus $vcfc(G)\leq2$, and
the proof is complete. \qed

From Theorem \ref{thm1} and Lemma \ref{lem2}, we have the following result.

\begin{cor}\label{cor1}
For the complete graph $K_n$ with $n\geq 2$, $vcfc(K_n)=2$.
\end{cor}

After the above preparation, graphs with $vcfc(G)=2$ can be
characterized.
\begin{thm}\label{thm2}
Let $G$ be a connected graph of order at least $3$. Then $vcfc(G)=2$
if and only if $G$ is $2$-connected or $G$ has only one cut vertex.
\end{thm}

\pf Firstly, we prove its sufficiency. If $G$ is $2$-connected, then
it follows from Lemma \ref{lem2} that $vcfc(G)=2$. Now suppose that
$G$ has exactly one cut vertex, say $w$. Since $vcfc(G)\geq 2$, we
just need to show that $vcfc(G)\leq 2$. Define a $2$-coloring $c$ of
the vertices of $G$ by coloring the vertex $w$ with color $2$ and
all the other vertices with color $1$. Since $G$ has only one cut
vertex, it follows that $G$ consists of some blocks which have the
common vertex $w$. Next it remains to check that for any two
vertices $u$ and $v$ in $G$, there is a conflict-free path between
them. It is clearly true for the case that $w\in\{u,v\}$. Thus we may
assume that $w\in V(G)\backslash\{u,v\}$. If $u$ and $v$ are in the
same block, then the block must be nontrivial. From Lemma \ref{lem1}
and the coloring $c$ of $G$, we get that there is a conflict-free
path from $u$ to $v$ in the block. If $u$ and $v$ are in two
different blocks, then there is a $u$-$w$ path $P_1$ and a $v$-$w$
path $P_2$ in the two blocks, respectively. Clearly, the path
$uP_1wP_2v$ is the desired path.

Now, we show its necessity. Let $vcfc(G)=2$. By Lemma \ref{lem2}, it
remains to show that if $G$ is not $2$-connected, then $G$ has only
one cut vertex. Suppose that $G$ has at least two cut vertices. Let
$B_1$ and $B_2$ be two blocks in $G$ which only contain one cut
vertex, respectively. Moreover, denote by $v_1$ and $v_2$ the cut
vertices in $B_1$ and $B_2$, respectively. Note that for any two
vertices in the same block, all paths connecting them are in the
block. Thus, each block needs two colors. Let $u_1$ be the vertex
whose color is different from $v_1$ in $B_1$ and $u_2$ be the vertex
whose color is different from $v_2$ in $B_2$. Clearly, all paths
from $u_1$ to $u_2$ in $G$ must pass through the vertices $v_1$ and
$v_2$. However, the four vertices $u_1,v_1,u_2,v_2$ use two colors
twice. Thus there does not exist a conflict-free path between $u_1$
and $u_2$ in $G$, a contradiction. Hence $G$ has only one cut
vertex. \qed

The following corollary is immediate from Theorem \ref{thm2}.

\begin{cor}\label{cor2}
Let $G$ be a connected graph. Then $vcfc(G)\geq 3$ if and only if $G$ has at least two cut vertices.
\end{cor}

Next we give two sufficient conditions for a graph $G$ to have  $vcfc(G)=3$.
\begin{thm}\label{thm3}
Let $G$ be a connected graph of order $n$ with maximum degree
$\Delta(G)$. If $G$ has at least two cut-vertices and
$n-4\leq\Delta(G)\leq n-2$, then $vcfc(G)=3$.
\end{thm}

\begin{figure}[htbp]
\begin{center}
\includegraphics[scale=1.0]{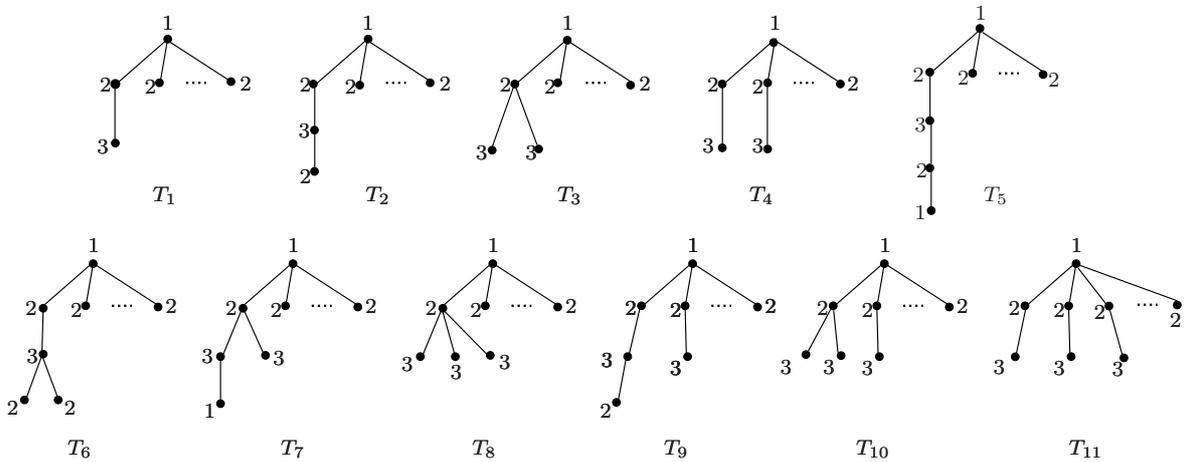}
\caption{The eleven graphs in Theorem 3.}\label{fig.1}
\end{center}
\end{figure}

\pf Since $G$ has at least two cut-vertices, it follows that
$vcfc(G)\geq 3$ by Corollary \ref{cor2}, and so we only need to show
that $vcfc(G)\leq 3$. We distinguish the following three cases to
show this theorem.

\textbf{Case 1.} $\Delta(G)=n-2$.

In this case, $G$ must have a spanning tree $T_1$ shown in Figure
\ref{fig.1}. Moreover, a $3$-coloring of the vertices of $T_1$ is
shown in Figure \ref{fig.1} to make $T_1$ conflict-free
vertex-connected. Thus $vcfc(T_1)\leq 3$. From Observation
\ref{obs1}, we have $vcfc(G)\leq vcfc(T_1)$, and hence $vcfc(G)\leq
3$.

\textbf{Case 2.} $\Delta(G)=n-3$.

Since $\Delta(G)=n-3$, it follows that $G$ must have a spanning tree
depicted as one of $T_i\ (2\leq i\leq4)$ shown in Figure
\ref{fig.1}. For $2\leq i\leq 4$, a $3$-coloring of the vertices of
$T_i$ is shown in Figure \ref{fig.1} to make $T_i$ conflict-free
vertex-connected. From Observation \ref{obs1}, we have $vcfc(G)\leq
vcfc(T_i)\leq3$.

\textbf{Case 3.} $\Delta(G)=n-4$.

Since $\Delta(G)=n-4$, it follows that $G$ must have a spanning tree
depicted as one of $T_i\ (5\leq i\leq11)$ shown in Figure
\ref{fig.1}. For $5\leq i\leq 11$, a $3$-coloring of the vertices of
$T_i$ is shown in Figure \ref{fig.1} to make $T_i$ conflict-free
vertex-connected. From Observation \ref{obs1}, we have $vcfc(G)\leq
vcfc(T_i)\leq 3$.

From the above argument, we conclude that $vcfc(G)=3$. \qed

\begin{figure}[htbp]
\begin{center}

\scalebox{0.6}[0.6]{\includegraphics{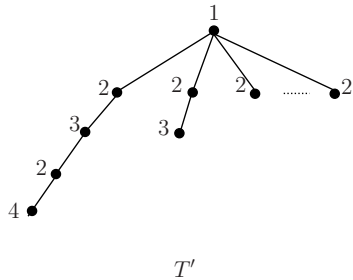}}
\end{center}
\caption{A $4$-coloring of the vertices of $T'$.}\label{fig.2}
\end{figure}

\begin{remark} The condition on the maximum degree above can not be improved, since if $G$ is $T'$ shown in Figure \ref{fig.2}, then $\Delta(G)=n-5$ and $vcfc(G)=4$. Note that there is only one path between any two vertices in a tree. Then any two adjacent vertices in $T'$ need two different colors. Considering this, we can check that three colors can not make $T'$ conflict-free vertex-connected and so $vcfc(T')\geq 4$. Moreover, a $4$-coloring of the vertices of $T'$ is shown in Figure \ref{fig.2} to make $T'$ conflict-free vertex-connected. Hence $vcfc(T')=4$.
\end{remark}

Let $C(G)$ denote the subgraph of $G$ induced by the set of cut-edges of $G$.

\begin{thm}\label{thm4}
Let $G$ be a connected graph with at least two cut-vertices. If
$C(G)$ is a star, then $vcfc(G)=3$.
\end{thm}

\pf By Corollary \ref{cor2}, it suffices to show that $vcfc(G)\leq
3$, since $G$ has at least two cut-vertices. Let
$V(C(G))=\{v_0,v_1,...,v_t\}$, where $t\geq1$ and $v_0$ is the
center of the star $C(G)$. Define a $3$-coloring of the vertices of
$G$ by coloring the vertex $v_0$ with color $1$, the pendant
vertices $\{v_1,...,v_t\}$ of $C(G)$ with color $2$ and all the
other vertices with color $3$. Next, it remains to check that for any
two vertices $u$ and $v$ in $G$, there is a conflict-free path
between them. If $u,v\in V(C(G))$, then the desired path is the
unique path from $u$ to $v$ in $C(G)$. If $u$ and $v$ belong to the
same nontrivial block, then by Lemma \ref{lem1}, there is a $u$-$v$
path in the block containing the vertex which is also in $C(G)$.
Clearly, this path is the desired path. Now we may assume that $u$
and $v$ are in two different nontrivial blocks. Consider a shortest
$u$-$v$ path in $G$. This path must go through the center $v_0$
which has the unique color $1$ and so it is the desired path. Thus,
$vcfc(G)\leq3$, and the proof is complete. \qed

The {\it $t$-corona} of a graph $H$, denoted by $Cor_t(H)$, is a
graph obtained from $H$ by adding $t$ pendant edges to each vertex
of $H$.
\begin{pro}\label{pro1}
Let $C_n$ be a cycle and $G$ be its $t$-corona, where $t\geq 1$. Then $vcfc(G)=3$.
\end{pro}

\pf Since $G$ has at least three cut-vertices, we have $vcfc(G)\geq
3$ by Theorem \ref{thm2}, and so it remains to show that
$vcfc(G)\leq 3$. Define a $3$-coloring $c$ of the vertices of $G$ by
coloring all the pendant vertices with color $1$, one of the
vertices of $C_n$ with color $2$ and the other vertices with color
$3$. It is easy to check that for any two vertices of $G$, there is
a conflict-free path between them. Then $vcfc(G)\leq 3$, and we
complete the proof. \qed

It seems that it is not easy to characterize graphs with
$vcfc(G)=3$. But, below we study a family of graphs with
conflict-free vertex-connection number three. Before it, we provide
the concept of a segment: Let $G$ be a connected graph whose block
graph is a path. Let $B_1, B_2,\cdots, B_k$ be the blocks of $G$
such that $|V(B_i)\bigcap V(B_{i+1})|=1$ and $E(B_i)\bigcap
E(B_{i+1})=\emptyset\ (1\leq i\leq k-1)$. We call $F_{p,q}\ (1\leq p\leq
q\leq k)$ a {\it segment} of $G$ if $F_{p,q}$ can be obtained from
$B_p,B_{p+1},\cdots,B_q$.

\begin{thm}\label{thm5}
Let $G$ be a connected graph with at least two cut-vertices, and its
block graph $B(G)$ is a path. Then $vcfc(G)=3$ if and only if $G$ is
a segment of one type of the thirteen graphs listed below.
 \end{thm}

\begin{figure}[htbp]
\begin{center}
\includegraphics[scale=0.9]{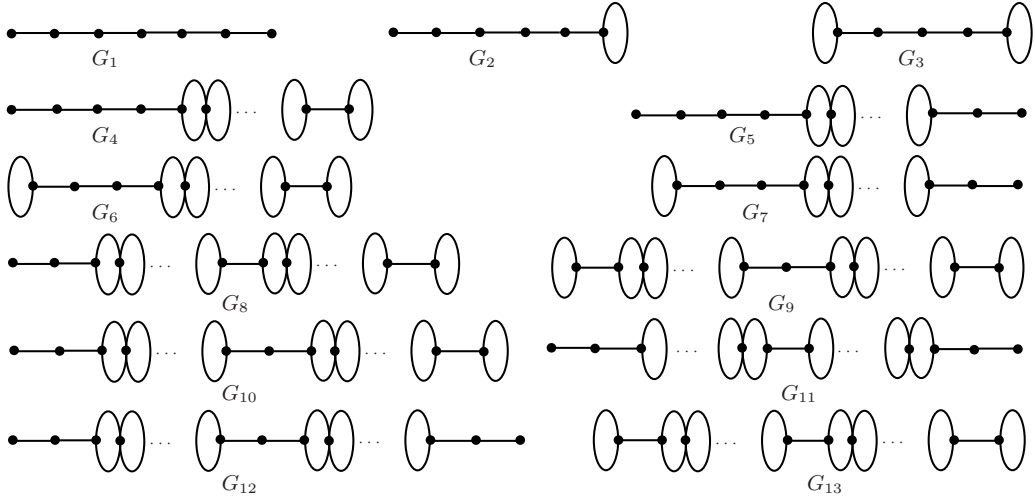}
\caption{The thirteen types of graphs}\label{fig3}
\end{center}
\end{figure}

\pf We divide the proof of this theorem into several cases in terms
of the maximum consecutive trivial blocks (cut-edges) in the graph $G$.
It can be easily found that $G$ can contain at most six
consecutive cut-edges since seven edges would require at least four
colors according to Theorem \ref{thm1}. By the way, we label the
vertices of $P_k$ as $v_1,v_2\cdots,v_k$ from left to right. And
$c(v)$ the color of the vertex $v$. We distinguish the following
cases to prove this theorem.

\textbf{Case 1:} $G$ contains a $P_7$.

Since $G$ contains a $P_7$, it follows that the unique color, say $3$,
must appear on $v_4$. If it appears on another position, it will
leave a path with at least four vertices colored only with two
colors. We can easily obtain that either of these two colors must
appear more than once in this path. Thus, no unique color is for this
path, a contradiction. Since $3$ appears on $v_4$ of $P_7$, it
leaves three vertices on the left and right sides, respectively. Since  $vcfc(P_3)=2$, the two $P_3$s both need colors $1$ and
$2$, and thus $3$ is the only unique color on $P_7$. If we add one
nontrivial block to an end of $P_7$, say $v_7$, then the block must
contain no vertex of color $3$ since the path between the
$3$-colored vertex in the block and $v_1$ would have both $1,2,3$
colors appear more than once, implying that there is no unique color for this path, a contradiction. But as a block must contain more than one colors, if we select a vertex $v$ in the nontrivial block such that $c(v)\neq c(v_7)$, then the path between $v$ and $v_5$ must contains no vertex of the unique color. In this way, $G$ can only be $P_7$ since all
other graphs would have a segment which consists of a $P_7$ with one
block to its end. Thus $G_1$ is the only graph for this case.

\textbf{Case 2:} $G$ contains a $P_6$ as its maximum consecutive
trivial blocks.

We know that $P_6$ surely requires three colors. The unique color of
$P_6$, say the color $3$, must appear on $v_3$ or $v_4$ for the same
reason as in Case 1. By symmetry, we assume that the color $3$ appears on $v_3$. It can be easily verified that $3$ is the only unique color
on $P_6$. Again, if we want to add one nontrivial block $B_1$ to an end of $P_6$, the block can not contain a vertex of color $3$. But
if this block is added to $v_6$, again if we choose a vertex $x$ of
the block such that $c(v)\neq c(v_6)$, the path between $x$ and
$v_4$ has no unique color. Thus the nontrivial block can be only
tied to the $v_1$ side if possible. Since $c(v_1)\neq 3$, we assume
that $c(v_1)=1$ without loss of generality. Then $c(v_2)=2$ and we
assign the color $2$ to every vertex of $B_1$ except for $v_1$. We
can easily verify that this graph is conflict-free vertex connected.
But if we add one more block after this nontrivial block $B_2$,
either trivial or nontrivial, both $B_1$ and $B_2$ must contain only
colors $1$ and $2$. Let the cut-vertex between $B_1$ and $B_2$ be
$y$, then pick a vertex $z$ of $B_2$ other than $y$ such that
$c(z)\neq c(y)$. We can check that the path between $z$ and $v_2$
has no unique color. Thus $G$ can only be $P_6$ or $P_6$ with one
nontrivial block tied to its one end as $G_2$ in Figure \ref{fig3}.

\textbf{Case 3:} $G$ contains a $P_5$ as its maximum consecutive
trivial blocks.

If $G$ contains a $P_5$ as its maximum consecutive trivial blocks,
then $P_5$ needs three colors according to Theorem \ref{thm1}. If
$3$ is the only unique color on $P_5$ and $c(v_3)=3$, then similar
to the case above, any block can not contain a vertex of color $3$
if later added. So we can not add more than one blocks to any end of
$P_5$. If we add two blocks $B_1$ and $B_2$ to one end of $P_5$, say
$v_5(v_5\in B_1)$, and the cut-vertex between $B_1$ and $B_2$ is
$x$. Then we can select a vertex $y$ from $B_2$ such that $c(y)\neq
c(x)$ and the path from $y$ to $v_4$ has no unique color. So we tie
one nontrivial block to $v_1$ and $v_5$ called $B_1$ and $B_5$,
respectively. Let $c(v_1)=c(v_5)=1$ and $c(v_2)=c(v_4)=2$, and
assign the color $2$ to any vertex in $B_1$ and $B_5$ except for
$v_1$ and $v_5$. Clearly, this is a conflict-free vertex-coloring of
$G$. And this corresponds to $G_3$ in Figure \ref{fig3}.

If $2$ and $3$ are two unique colors on $P_5$, then without loss of
generality, we set $c(v_2)=2$, $c(v_4)=3$ and
$c(v_1)=c(v_3)=c(v_5)=1$. Then if the nontrivial blocks added to one
side of $P_5$ contain both colors $2$ and $3$, the other side can
not be tied to any nontrivial block containing the color $2$ or $3$.
Otherwise, a path with no unique color will be found. But this is
impossible since any block must contain at least two colors. We
conclude that if we add nontrivial blocks to both sides of $P_5$,
then these nontrivial blocks can only contain the color $1$ and only
one color of colors $2$ and $3$. We assume that they only contain the
colors $1$ and $2$, without loss of generality. But then let $B$ be
the nontrivial block containing $v_1$ and $x\in V(B)$ with $c(x)=2$,
consequently the path between $x$ and $v_3$ has no unique color, a
contradiction. Thus we conclude that we can only add nontrivial
blocks to one end of $P_5$, say $v_5$, and we will see the number is
not limited. So we can now label these nontrivial blocks from the
side near $v_5$ as $B_1,B_2,\cdots,B_k$. Pick one vertex from $B_i$
other than a cut-vertex and give it the color $2$ if $i$ is odd.
Then pick one vertex other than a cut-vertex from $B_i$ and assign
the color $3$ to it if $i$ is even. Finally, assign the color $1$ to
all other unmentioned vertices in $B_i \ (1\leq i\leq k)$. We can
easily check that this coloring is conflict-free vertex-connected.
Remember that we can not add two segments of trivial blocks with
nontrivial blocks between them after $B_k$, since two segments of
trivial blocks would both contain the color $1$ and only one color
of $2$ and $3$, say the color $2$. Then the nontrivial blocks
between them must contain one vertex $x$ such that $c(x)=3$.
Therefore, the path between $x$ and $v_1$ would have no unique color,
a contradiction. By the way, we can not add $P_j \ (j\geq4)$ after
$B_k$ since $P_j$ would require more than two colors and thus the path
containing both $P_j$ and $P_5$ would have no unique color.
Therefore, we can only add $P_3$ or $P_2$ after $B_k$.

If we add $P_3$ after $B_k$ and label the vertices from the side
near $B_k$ as $u_1$, $u_2$ and $u_3$. Then $P_3$ must contain the
color $1$ and only one of color $2$ and $3$, say $2$. Then clearly
we can not add one more nontrivial block $B$ after $u_3$ since $B$
does not contain a vertex of color $3$ according to the discussion
above. Select a vertex $z$ of $B$ such that $c(z)\neq c(u_3)$ and
the path between $z$ and $u_1$ will have no unique color. Then if
$k$ is odd, set $c(u_1)=c(u_3)=1$, $c(u_2)=3$; and
$c(u_1)=c(u_3)=1$, $c(u_2)=2$ otherwise. The coloring of other
vertices is as the description above. It can be verified that this
coloring is conflict-free vertex-connected.

If we add $P_2$ after $B_k$ and label the vertices from the side
near $B_k$ as $u_1$ and $u_2$. $P_2$ must occupy two colors, say $1$
and $2$, then we can not add more than one nontrivial block after
$u_2$. Since if we add two nontrivial blocks $A_1$ and $A_2$ with
$x$ being the cut-vertex between them such that $u_2\in A_1$, bear
in mind that $A_1$ and $A_2$ can only contain colors $1$ and $2$.
Choose one vertex $y$ other than $x$ from $A_2$ such that $c(y)\neq
c(x)$. Then the path between $y$ and $u_1$ has no unique color, a
contradiction. Now we add one nontrivial block $A$ containing $u_2$.
So we let $c(u_2)=2$ if $k$ is even and $c(u_2)=3$ otherwise.
Besides, let $c(u_1)=1$ and assign the color $1$ to all other
vertices in $A$ except for $u_2$. We can check that this coloring is
conflict-free vertex-connected. The corresponding graphs are $G_4$
and $G_5$.

For the remaining case, when $c(v_4)=3$, $c(v_1)=c(v_3)=1$ and
$c(v_2)=c(v_5)=2$. We can check that the graph structure that allows
$3$ colors to make it conflict-free vertex-connected, has been covered by the discussion above. Thus we finish the proof of Case $3$.

Other cases can be similarly dealt with to the discussion above. But
as the process is rather complicated. We omit it here.\qed

At the end of this section, we pose the following problem.

\begin{prob}\label{prob2}
Characterize all the graphs $G$ with $vcfc(G)=3$.
\end{prob}

\section{Trees}

A {\it $k$-ranking} of a connected graph $G$ is a labeling of its
vertices with labels $1,2,3,\cdots,k$ such that every path between
any two vertices with the same label $i$ in $G$ contains at least
one vertex with label $j>i$. A graph $G$ is said to be {\it
$k$-rankable} if it has a $k$-ranking. The minimum $k$ for which $G$
is $k$-rankable is denoted by $r(G)$.

Iyer \cite{IRV} obtained the following result.
\begin{lem}{\upshape \cite{IRV}}\label{lem3}
Let $T$ be a tree of order $n\geq 3$. Then $r(T)\leq \log_{\frac{3}{2}}n$.
\end{lem}

The next two lemmas are preparations for Theorem \ref{thm6}.
\begin{lem}\label{lem4}
Let $G$ be a connected graph. Then $vcfc(G)\leq r(G)$.
\end{lem}

\pf Consider a ranking of the vertices of $G$. For any two vertices
$u$ and $v$ of $G$, let $P$ be a path between them and $k$ be the
maximum label of the vertices of $P$. If there is only one vertex with
label $k$ in $P$, then the proof is done. So we assume that $P$
contains at least two vertices with label $k$. According to the
definition of ranking, there must exist one vertex with label $j>k$
on $P$, which is a contradiction. Hence $P$ contains only one vertex
with label $k$. View the $r(G)$-ranking of $G$ as its
vertex-coloring with $r(G)$ colors. Then the path $P$ is a
conflict-free path between $u$ and $v$ in $G$. Thus $vcfc(G)\leq
r(G)$.\qed

\begin{lem}\label{lem5} Let $T$ be a nontrivial tree. Then $vcfc(T)\geq \chi(T)$, where $\chi(T)$ is the chromatic number of $T$ and the bound is sharp.
\end{lem}

\pf Define a vertex-coloring of $T$ with $vcfc(T)$ colors such that
$T$ is conflict-free vertex-connected. Since there is only one path
between any two vertices in $T$, it follows that any two adjacent
vertices must have different colors, and hence $vcfc(T)\geq \chi(T)$.
To show the sharpness of the bound, we let
$T$ be a star of order at least two. Clearly, $\chi(T)=2$. By Theorem
\ref{thm2}, we have $vcfc(T)=2(=\chi(T))$. \qed

Combining the lemmas above, we can have the following bounds for $vcfc(T)$ of a tree $T$.

\begin{thm}\label{thm6}
Let $T$ be a tree of order $n\geq 3$ and $d(T)$ be its diameter. Then
\begin{eqnarray*}
\max\{\chi(T),\lceil\log_{2}(d(T)+1)\rceil\}\leq vcfc(T)\leq
\log_{\frac{3}{2}}n.
\end{eqnarray*}
\end{thm}

\pf The lower bound is an immediate result from Lemma
\ref{lem5} and Theorem \ref{thm1}, while the upper bound can be
deduced from Lemmas \ref{lem3} and \ref{lem4}.\qed

Let $G$ be a connected graph. The {\it eccentricity} $\epsilon_G(v)$
of a vertex $v$ in $G$ is the maximum value among the distances between $v$ and the other vertices in $G$. The {\it radius} $rad(G)$ of $G$ is
the minimum eccentricity among all the vertices of $G$. A {\it central vertex} of radius $rad(G)$ is one whose eccentricity is $rad(G)$. Remind that $d_G(u,v)$ is the shortest distance between the two vertices $u$ and $v$ in $G$.

\begin{thm}\label{thm7}
Let $T$ be a tree with radius $rad(T)$. Then $vcfc(T)\leq rad(T)+1$.
Moreover, the bound is sharp.
\end{thm}

\pf Let $v$ be a central vertex of radius $rad(T)$ in $T$. Let
$V_i=\{u\in V(T):\ d_T(u,v)=i\}$, where $0\leq i\leq rad(T)$. Hence
$V_0=\{v\}$. Define a vertex-coloring $c$ of $T$ with $rad(T)+1$
colors by coloring the vertices of $V_i$ with color $i+1$, where
$0\leq i\leq rad(T)$. It is easy to check that for any two vertices
of $T$, there is a conflict-free path between them, and hence
$vcfc(T)\leq rad(T)+1$. To show the sharpness of the bound, we let
$T$ be a star of order at least two. Clearly, $rad(T)=1$. By Theorem
\ref{thm2}, we have $vcfc(T)=2(=rad(T)+1)$. \qed

For each connected graph $G$, we can always find a spanning tree $T$
of $G$ such that $rad(T)=rad(G)$. From Observation \ref{obs1} and
Theorem \ref{thm7}, we can get the following result.

\begin{cor}\label{cor3}
Let $G$ be a connected graph. Then $vcfc(G)\leq rad(G)+1$.
\end{cor}

For trees, we can give an upper bound of its conflict-free
vertex-connection number in term of its order.
\begin{pro}\label{pro2}
Let $T$ be a tree with order $n\geq 5$. Then
$vcfc(T)\leq\lceil\frac{n}{2}\rceil$. Moreover, the bound is sharp.
\end{pro}

\pf If $T$ is a path, then it follows from Theorem \ref{thm1} that
$vcfc(T)=\lceil\log_2(n+1)\rceil\leq\lceil\frac{n}{2}\rceil$. From
now on, we suppose that $T$ is not a path. Then the longest path in
$T$ has at most $n-1$ vertices. So we have $rad(T)\leq\frac{n-1}{2}$
if $n$ is odd and $rad(T)\leq\frac{n-2}{2}$ if $n$ is even. By
Theorem \ref{thm7}, we have $vcfc(T)\leq rad(T)+1$, and hence
$vcfc(T)\leq\lceil\frac{n}{2}\rceil$. To show the sharpness of the
bound, we set $T=P_5$. Then $vcfc(T)=3$ by Theorem \ref{thm1} and
$\lceil\frac{n}{2}\rceil=3$. \qed

To end this section, we give the following conjecture.

\begin{conj}\label{conj1}
Let $G$ be a connected graph of order $n$. Then $vcfc(G)\leq vcfc(P_n)$.
\end{conj}

\end{document}